\documentclass[10pt,a4paper,english]{article}
\usepackage[T1]{fontenc}
\usepackage[latin9]{inputenc}
\usepackage{amsthm}
\usepackage{amsmath}
\usepackage{amssymb}
\usepackage[authoryear]{natbib}

\theoremstyle{plain}
\theoremstyle{plain}
\newtheorem{thm}{Theorem}
  \theoremstyle{plain}
  \newtheorem{lem}[thm]{Lemma}
  \theoremstyle{plain}
  \newtheorem{cor}[thm]{Corollary}

 \newtheorem{hypothesis}{Hypothesis}

\providecommand{\tv}[1]{\lVert #1 \rVert_{TV}}

\usepackage{babel}

\usepackage{babel}

\begin{document}

\title{Stability of Feynman-Kac formulae with path-dependent potentials}
\date{}

\author{N. Chopin%
\thanks{Corresponding author: nicolas.chopin@ensae.fr, ENSAE-CREST, 3, Avenue Pierre Larousse, 92240 Malakoff, France%
}, Pierre Del Moral%
\thanks{INRIA Bordeaux Sud-Ouest, 351, cours de la Lib\'eration
33405 Talence cedex, France%
} and S. Rubenthaler%
\thanks{Laboratoire J.-A. Dieudonn\'e, Universit\'e de Nice-Sophia Antipolis,
Parc Valrose, 06108 Nice, cedex 2, France %
}}
\maketitle
\begin{abstract}
Several particle algorithms admit a Feynman-Kac representation such
that the potential function may be expressed as a recursive function
which depends on the complete state trajectory. An important example
is the mixture Kalman filter, but other models and algorithms of practical
interest fall in this category. We study the asymptotic stability
of such particle algorithms as time goes to infinity. As a corollary,
practical conditions for the stability of the mixture Kalman filter,
and a mixture GARCH filter, are derived. Finally, we show that our
results can also lead to weaker conditions for the stability of standard
particle algorithms, such that the potential function depends on the
last state only.
\end{abstract}

\section{Introduction}

The most common application of the theory of Feynman-Kac formulae
\citep[see e.g.][]{delMoral:book} is nonlinear filtering of a hidden
Markov chain $(\Lambda_{n})$, based on observed process $(Y_{n})$.
In such settings, the potential function at time $n$ typically depends
only on the current state $\Lambda_{n}$. The uniform stability of
the corresponding particle approximations can be obtained under appropriate
conditions, see Section 7.4.3 of the aforementioned book and references
therein. For a good overview of the theoretical and methodological aspects
of particle approximation algorithms, also known as particle filtering
algorithms, see also \citet{DouFreiGor}, \citet{Kun:SSHMM}, 
and \citet{CapMouRyd}.

They are however several applications of practical interest where the
potential function depends on the complete state trajectory $\Lambda_{0:n}=(\Lambda_{0},\ldots,\Lambda_{n})$.
The corresponding particle filtering algorithms still have a fixed
computational cost per iteration, because the potential can be computed
using recursive formulae. An important example is the class of conditional
linear Gaussian dynamic models, where the conditioning is on some
unobserved Markov chain $\Lambda_{n}$. The corresponding particle
algorithm is known as the mixture Kalman filter (\citealp{ChenLiu},
see also Example 7 in \citealp{DouGodAnd}, and \citealp{AndDou}, for
a related algorithm): the potential function at time $n$ is then
a Gaussian density, the parameters of which are computed recursively
using the Kalman-Bucy filter \citep{KalBuc}. Another example is the
mixture GARCH model considered in \citet{Chopin:change}.

It is worth noting that these models such that the potential functions
are path-dependent can often be reformulated as a standard hidden
Markov model, with a potential function depending on the last state
only, by adding components to the hidden Markov chain. For instance,
the mixture Kalman filter may be interpreted as a standard particle
filtering algorithm, provided the hidden Markov process is augmented
with the associated Kalman filter parameters (filtering expectation
and error covariance matrix) that are computed iteratively in the
algorithm. However, this representation is unwieldy, and the augmented
Markov process does not fulfil the usual mixing conditions found in
the literature on the stability of particle approximations. This
is the main reason why our study is based on path-dependent potential
functions. Quite interestingly, we shall see that the opposite perspective
is more fruitful. Specifically, our stability results obtained for
path-dependent potential functions can also be applied to standard
state-space models, leading to stability results under 
conditions different from those previously given in the literature. 

In this paper, we study the asymptotic stability of particle algorithms
based on path-dependent potential functions. We work under the assumption
that the dependence of potential $n$ on state $n-p$ vanishes exponentially
in $p$. This assumption is met in practical settings because of the
recursive nature of the potential functions. Our proofs are based
on the following construction: the true filter is compared with an
approximate filter associated to `truncated' potentials, that is potentials
that depend only on $\lambda_{n-p+1:n}$, the vector of the last $p$
states, for some well-chosen integer $p$. Then, we compare the truncated
filter with its particle approximation, using the fact the `truncated'
filter corresponds to a standard Feynman-Kac model with a Markov chain
of fixed dimension. Finally, we use a coupling construction to compare
the particle approximations of the true filter and the truncated filter.
In this way, we obtain estimates of the stability of the particle
algorithm of interest. We apply our results to the two aforementioned
classes of models, and obtain practical conditions under which the
corresponding particle algorithms are stable uniformly in time.

The paper is organised as follows. Section \ref{sec:model} introduces
the model and the notations. Section \ref{sec:local} evaluates the
local error induced by the truncation. Section \ref{sec:mix-truncated}
studies the mixing properties of the truncated filter. Section \ref{sec:Propagation}
studies the propagation of the truncation error. Section \ref{sec:Coupling}
develops a coupling argument for the two particle systems. Section
\ref{sec:Main} states the main theorem of the paper, which provides
a bound for the particle error and  derives time-uniform
estimates for the long-term propagation of the error in the particle
approximation of the true model. Section \ref{sec:applis} applies these
results to two particle algorithms of practical interest, namely,
the mixture Kalman filter, and the mixture GARCH filter, and shows
how these results can be adapted to standard state-space models, such
that the potential function depends only on the last state.

\section{Model and notations\label{sec:model}}

We consider a hidden Markov model, with latent (non-observed) state
process $\{\Lambda_{n},n\geq0\}$, and observed process
$\{Y_{n},n\geq1\}$, taking values respectively in a complete separable
metric space $E$ and in $F=\mathbb{R}^{d}$. The state process is an
inhomogeneous Markov chain, with initial probability distribution
$\zeta$, and transition kernel $Q_{n}$. The observed process $Y_{n}$
admits $\Psi_{n}(y_{n}|y_{1:n-1};\lambda_{0:n})$ as a conditional
probability density (with respect to an appropriate dominating
measure) given $\Lambda_{0:n}=\lambda_{0:n}$ and
$Y_{1:n-1}=y_{1:n-1}$, where the short-hand $v_{0:n}$ for any symbol
$v$ stands for the vector $(v_0,\ldots,v_n)$.  As explained in the
Introduction, this quantity depends on the entire path
$\lambda_{0:n}$, rather than the last state $\lambda_{n}$. Following
common practice, we drop dependencies on the $y_{n}$'s in the
notations, as the observed sequence $y_{0:n}$ may be considered as
fixed, and use the short-hand
$\Psi_{n}(\lambda_{0:n})=\Psi_{n}(y_{n}|y_{0:n-1};\lambda_{0:n})$.
The model admits a Feynman-Kac representation which we describe fully
in (\ref{Eq:FKrep}). We consider the following assumptions.

\begin{hypothesis} For all $n\geq 1$, the kernel $Q_{n}$ is mixing, i.e.
there exists $\varepsilon_{n}\in(0,1)$ such that
 \[
\varepsilon_{n}\xi(A)\leq Q_{n}(\lambda_{n-1},A)\leq\frac{1}{\varepsilon_{n}}\xi(A)\]
 for some $\xi\in\mathcal{M}_{+}(E)$, and for any Borel set $A\subset E$,
any $\lambda_{n-1}\in E$. 
\end{hypothesis} 

\begin{hypothesis} 
For $p$ large enough, and all $n\geq p$, there exists a `truncated'
potential function $\tilde{\Psi}_{n}^{p}(\lambda_{n-p+1:n})$ that
depends on the last $p$ states only, and that approximates $\Psi_{n}$
in the sense that
\[
|\Psi_{n}(\lambda_{0:n})-\tilde{\Psi}_{n}^{p}(\lambda_{n-p+1:n})|
\leq\phi_{n}\tau^{p}\left\{ \Psi_{n}(\lambda_{0:n})\wedge\tilde{\Psi}_{n}^{p}(\lambda_{n-p+1:n})\right\} \]
 for some constants $\phi_{n}$ and $\tau$, $\phi_{n}>0$, $0<\tau<1$,
and all $\lambda_{0:n}\in E^{n+1}$. For convenience, we abuse notations
and set $\tilde{\Psi}_{n}^{p}=\Psi_{n}$ for $p>n$.
\end{hypothesis} 

\begin{hypothesis} 
There exists constants $a_{n},b_{n}$, $n\geq0$,
$a_{n}\geq1$, $b_{n}\geq1$, such that
 \[
\frac{1}{a_{n}}\leq\Psi_{n}(\lambda_{0:n})\leq b_{n}
,\quad\frac{1}{a_{n}}\leq\tilde{\Psi}_{n}^{p}(\lambda_{(n-p+1)^{+}:n})\leq b_{n}\]
 for all $\lambda_{0:n}\in E^{n+1}$, using the short-hand $k^{+}=k\vee 0$
for any integer $k$.

\end{hypothesis}

The constants $a_{n}$ and $\phi_{n}$  depend implicitly on the
realisation $y_{1:n}$ of the observed process. Hypotheses 1 and 3
are standard in the filtering literature; see e.g. \citet[][]{delMoral:book}. Hypothesis 2 formalises
the fact that potential functions are computed using iterative formulae,
and therefore should forget past states at an exponential rate. One
may take $\tilde{\Psi}_{n}^{p}(\lambda_{n-p+1:n})=\Psi_{n}(z,\ldots,z,\lambda_{n-p+1:n})$
for instance, where $z$ is an arbitrary element of $E$. We shall
work out, in several models of interest, practical conditions under which Hypothesis
2 is fulfilled in Section \ref{sec:applis}.

We introduce the following notations for the forward kernels, for
$n\geq1$:
 
$$  \gamma_{n}(\lambda_{0:n-1},d\lambda_{0:n}') = 
 \delta_{\lambda_{0:n-1}}(d\lambda_{0:n-1}')Q_n(\lambda_{n-1},d\lambda'_{n})
\Psi_{n}(\lambda'_{0:n})
$$
where $\delta_{\lambda_{0:n-1}}$ is the
Dirac measure centred at $\lambda_{0:n-1}$.  The above kernels
implicitly defines operators on measures and on test functions,
i.e.,
$$
 \gamma_{n}\mu(f)=\left\langle \gamma_{n}\mu,f\right\rangle = 
    \int\mu(d\lambda_{0:n-1})\gamma_n(\lambda_{0:n-1},d\lambda_{0:n}')  f(\lambda_{0:n}'),
$$
for any $\mu\in\mathcal{M}_{+}(E^{n+1})$, any
test function $f:E^{n+1}\rightarrow[0,1]$,
where $\mathcal{M}_{+}(E^{k})$ denotes the set of nonnegative
measures w.r.t. $E^{k}$, and $\mathcal{P}(E^{k})$ the set of probability
measures w.r.t. $E^{k}$. 

We associate to $\gamma_n$ a ``normalised'' operator $R_n$, such that,
for any $\mu\in  \mathcal{M}_{+}(E^n)$, $R_n\mu$ is defined as: 
$$ R_n\mu(f) = \frac{\gamma_n\mu(f)}{\gamma_n\mu(1)}$$
for any $f:E^{n+1}\rightarrow\mathbb{R}^{+}$.  Both the $\gamma_n$'s
and the $R_n$'s may be iterated using the following short-hands, for $1\leq
k\leq n$:
$$ \gamma_{k:n}\mu=\gamma_{n}\ldots\gamma_{k}\mu,\quad
R_{k:n}\mu=R_{n}\ldots R_{k}\mu.$$

We have the following Feynman-Kac representation:
\begin{equation}
\mathbb{E}(f(\Lambda_{0:n})|Y_{1:n}=y_{1:n})=R_{1:n}\zeta(f)\ ,
\label{Eq:FKrep}
\end{equation}
$\forall n$, $\forall f:E^{n+1}\rightarrow\mathbb{R}^{+}$, where,
as mentioned above, $\zeta$ the law of $\Lambda_{0}$.

Finally, we denote the total variation norm on nonnegative measures by
$\tv{\cdot}$, the supremum norm on bounded functions by
$\Vert\cdot\Vert_{\infty}$, and the Hilbert metric by $h(\mu,\mu')$
for any pair $\mu,\mu'\in\mathcal{M}_{+}(E^{k})$, $k\geq1$; see
e.g. \citet{atar1997esn} or \citet{legland2004sau}, Definition 3.3. We
recall that the Hilbert metric is scale invariant, and is related to
the total variation norm in the following way, see e.g. Lemma 3.4 in
\citet{legland2004sau}:
\begin{eqnarray}
  \tv{\mu-\mu'} & \leq & \frac{2}{\log3}h(\mu,\mu')\label{Eq:borne1}\\
  h(K\mu,K\mu') & \leq & \frac{1}{\varepsilon^{2}}\tv{\mu-\mu'}\label{Eq:borne2}
\end{eqnarray}
provided $K$ is a $\varepsilon$-mixing kernel. We can also derive
the following properties from the definition of $h$ ($\forall k\in\mathbb{N}^{*}$,
$\forall\mu,\mu'\in\mathcal{M}(E^{k})$): 
\begin{eqnarray}
  &  & \text{ \ensuremath{\forall}kernel \ensuremath{Q}, }
  h(Q\mu,Q\mu')\leq h(\mu,\mu')\ ,\label{Eq:contract1}\\
  &  & \text{\ensuremath{\forall}nonnegative function
    \ensuremath{\psi}, }h(\psi\mu,\psi\mu')\leq h(\mu,\mu')\label{Eq:inv1}
\end{eqnarray}
with an equality in the latter equation if $\psi$ is positive.

\section{Local error induced by truncation\label{sec:local}}

Until further notice, $p$ is a fixed integer such that $p\geq2$
and such that Hypothesis 2 holds. Since our proofs involve a comparison
between the true filter and a `truncated' filter, we introduce the
projection operator $H_{n}^{p}$ which, for $n\geq p$, associates
to any measure $\mu(d\lambda_{0:n})\in\mathcal{M}_{+}(E^{n+1})$ its
marginal w.r.t. its last $p$ components, i.e. :
 \[
H_{n}^{p}(\mu)(f)=\int\mu(d\lambda_{0:n})f(\lambda_{n-p+1:n})\]
 for any $f:E^{p}\rightarrow\mathbb{R}$; for $p>n$, let $H_{n}^{p}(\mu)=\mu$.
We also define the following `truncated' forward kernels, for $ n\geq p$:
\begin{eqnarray*}
 \lefteqn{\tilde{\gamma}_{n}^{p}(\lambda_{n-p:n-1},d\lambda_{n-p+1:n}')} \\
& = & 
\delta_{\lambda_{n-p+1:n-1}}(d\lambda_{n-p+1:n-1}')Q_{n}(\lambda_{n-1},d\lambda_{n}')
\tilde{\Psi}_{n}^{p}(\lambda'_{n-p+1:n})
\end{eqnarray*}
 and the associated  normalised operators, for $\mu\in\mathcal{M}_{+}(E^p)$,
 $f:E^p\rightarrow\mathbb{R}^{+}$:
 \[
\tilde{R}_{n}^{p}\mu(f) = 
\frac{\tilde{\gamma}_{n}^{p}\mu(f)}
{\tilde{\gamma}_{n}^{p}\mu(1)}
\]
 and set $\tilde{\gamma}_{n}^{p}=\gamma_{n}$, $\tilde{R}_{n}^{p}=R_{n}$
for $n<p$. From now on, we will refer to the filter associated to
these `truncated' operators as the truncated filter. 

We now evaluate the local error induced by the truncation.
\begin{lem}
\label{lem:loc2} For all $1\leq k<n$, and for all $\mu\in\mathcal{M}_{+}(E^{k})$,
\[
\left\Vert \tilde{R}_{k+1:n}^{p}H_{k}^{p}R_{k}\mu-
\tilde{R}_{k:n}^{p}H_{k-1}^{p}\mu\right\Vert _{TV}\leq2\phi_{k}\tau^{p}.\]
\end{lem}
\begin{proof}
Let $f:E^{p\wedge(n+1)}\rightarrow[0,1]$. One has
 \begin{eqnarray*}
\tilde{R}_{k+1:n}^{p}H_{k}^{p}R_{k}\mu(f) 
& = & \frac{\tilde{\gamma}_{k+1:n}^{p}H_{k}^{p}\gamma_{k}\mu(f)}{\tilde{\gamma}_{k+1:n}^{p}H_{k}^{p}\gamma_{k}\mu(1)}\\
\tilde{R}_{k:n}^{p}H_{k-1}^{p}\mu(f) 
& = & \frac{\tilde{\gamma}_{k:n}^{p}H_{k-1}^{p}\mu(f)}{\tilde{\gamma}_{k:n}^{p}H_{k-1}^{p}\mu(1)}
\end{eqnarray*}
where 
\begin{eqnarray*}
\tilde{\gamma}_{k+1:n}^p H_k^{p}\gamma_{k}\mu(f) 
& = & \int_{E^{n+1}}\mu(d\lambda_{0:k-1})Q_{k}(\lambda_{k-1},d\lambda_{k})\Psi_{k}(\lambda_{0:k})f(\lambda_{(n-p+1)^{+}:n})\\
 &  & \quad\times\prod_{i=k+1}^{n}\left[Q_{i}(\lambda_{i-1},d\lambda_{i})\tilde{\Psi}_{i}^{p}(\lambda_{(i-p+1)^{+}:i})\right]
\end{eqnarray*}
and 
\begin{eqnarray*}
\tilde{\gamma}_{k:n}^p H_{k-1}^{p}\mu(f) 
& = & \int_{E^{n+1}}\mu(d\lambda_{0:k-1})Q_{k}(\lambda_{k-1},d\lambda_{k})\tilde{\Psi}_{k}^{p}(\lambda_{k-p+1:k})f(\lambda_{(n-p+1)^{+}:n})\\
 &  & \quad\times\prod_{i=k+1}^{n}\left[Q_{i}(\lambda_{i-1},d\lambda_{i})\tilde{\Psi}_{i}^{p}(\lambda_{(i-p+1)^{+}:i})\right]
\end{eqnarray*}
hence 
\begin{eqnarray*}
\lefteqn{\left|\tilde{\gamma}_{k+1:n}H_{k}^{p}\gamma_{k}\mu(f)-\tilde{\gamma}_{k:n}H_{k-1}^{p}\mu(f)\right|}\\
& \leq & \int_{E^{n+1}}\mu(d\lambda_{0:k-1})Q_{k}(\lambda_{k-1},d\lambda_{k})\left|\Psi_{k}(\lambda_{0:k})-\tilde{\Psi}_{k}^{p}(\lambda_{(k-p+1)^{+}:k})\right|\\
 &  & \qquad\times f(\lambda_{(k-p+1)^{+}:k})\prod_{i=k+1}^{n}\left[Q_{i}(\lambda_{i-1},d\lambda_{j})\tilde{\Psi}_{i}^{p}(\lambda_{(i-p+1)^{+}:i})\right]\\
 & \leq & \phi_{k}\tau^{p}\int_{E^{n+1}}\mu(d\lambda_{0:k-1})Q_{k}(\lambda_{k-1},d\lambda_{k})(\Psi_k(\lambda_{0:k})\wedge\tilde{\Psi}_{k}^{p}(\lambda_{(k-p+1)^{+}:k}))\\
 &  & \quad\times f(\lambda_{(k-p+1)^{+}:k})\prod_{i=k+1}^{n}\left[Q_{i}(\lambda_{i-1},d\lambda_{i})\tilde{\Psi}_{i}^{p}(\lambda_{(i-p+1)^{+}:i})\right]\\
 & \leq & \phi_{k}\tau^{p}\left\{ \tilde{\gamma}_{k+1:n}H_{k}^{p}\gamma_{k}\mu(f)\wedge\tilde{\gamma}_{k:n}H_{k-1}^{p}\gamma_{k-1}\mu(f)\right\} 
\end{eqnarray*}
according to Hypothesis 2. And, since, for all $a$, $b$, $c$,
$d\in\mathbb{R}^{+}$ such that $a/b\leq1$ and $c/d\leq1$,

\begin{eqnarray}
\left|\frac{a}{b}-\frac{c}{d}\right| & \leq & \frac{|a-c|}{b}+\frac{|d-b|}{b}\label{eq:abcd}\end{eqnarray}
 one may conclude directly by taking $a=\tilde{\gamma}_{k+1:n}H_{k}^{p}\gamma_{k}\mu(f)$,
$b=\tilde{\gamma}_{k:n}H_{k-1}^{p}\gamma_{k-1}\mu(1)$, $c=\tilde{\gamma}_{k+1:n}^{p}H_{k}^{p}\mu(f)$,
and $d=\tilde{\gamma}_{k:n}^{p}H_{k-1}^{p}\mu(1)$. \end{proof}
\begin{lem}
\label{lem:32} For $k\geq1$, if there exists a (possibly random)
probability kernel $\bar{R}_{k}:E^{k\wedge p}\rightarrow\mathcal{P}(E^{(k+1)\wedge p})$
such that, for all $\mu\in\mathcal{P}(E^{k\wedge p})$, \[
\sup_{f:\Vert f\Vert_{\infty}=1}\mathbb{E}\left(\left|\langle\tilde{R}_{k}^{p}\mu-\bar{R}_{k}\mu,f\rangle\right|\right)\leq\delta_{k}\]
 for some $\delta_{k}\geq0$, then, for all $i\geq1$ and $\mu\in\mathcal{P}(E^{k\wedge p})$,
 \[
\sup_{f:\Vert f\Vert_{\infty}=1}\mathbb{E}\left(\left|\langle\tilde{R}_{k:k+i}^{p}\mu-\tilde{R}_{k+1:k+i}^{p}\bar{R}_{k}\mu,f\rangle\right|\right)\leq2(a_{k+1}\dots a_{k+i})(b_{k+1}\dots b_{k+i})\delta_{k}\]
 where the expectation is with respect to the distribution of $\bar{R}_{k}$.\end{lem}
\begin{proof}
Using the same ideas as above, one has, for $f:E^{(k+1-p+1)\wedge
  p}\rightarrow[0,1]$, 
\[
\langle\tilde{R}_{k:k+i}^{p}\mu-\tilde{R}_{k+1:k+i}^{p}\bar{R}_{k}\mu,f\rangle=\frac{\tilde{\gamma}_{k+1:k+i}^{p}\tilde{R}_{k}^{p}\mu(f)}{\tilde{\gamma}_{k+1:k+i}^{p}\tilde{R}_{k}^{p}\mu(1)}-\frac{\tilde{\gamma}_{k+1:k+i}^{p}\bar{R}_{k}\mu(f)}{\tilde{\gamma}_{k+1:k+i}^{p}\bar{R}_{k}\mu(1)}.\]
In order to use inequality (\ref{eq:abcd}), compute
 \begin{eqnarray*}
\lefteqn{\mathbb{E}(\left|\tilde{\gamma}_{k+1:k+i}^{p}\tilde{R}_{k}^{p}\mu(f)-\tilde{\gamma}_{k+1:k+i}^{p}\bar{R}_{k}\mu(f)\right|)}\\
 & = & \mathbb{E}\left(\left|\int(\tilde{R}_{k}^{p}\mu-\bar{R}_{k}\mu)(d\lambda_{(k-p+1)^{+}:k})\right.\right.\\
 &  & ~~~~~\left.\left.\prod_{l=k+1}^{k+i}Q_{l}(\lambda_{l-1},d\lambda_{l})\tilde{\Psi}_{l}^{p}(\lambda_{(l-p+1)^{+}:l})f(\lambda_{(k+i-p+1)^{+}:k+i})\right|\right)\\
 & \leq & \mathbb{E}\left(b_{k+1}\ldots b_{k+i}\left|(\tilde{R}_{k}^{p}\mu-\bar{R}_{k}\mu)(\bar{f})\right|\right)\\
 & \leq & b_{k+1}\ldots b_{k+i}\delta_{k}
\end{eqnarray*}
 where $\bar{f}$ is defined as
 \[
\bar{f}(\lambda_{(k-p+1)^{+}:k})=\int\prod_{l=k+1}^{k+i}Q_{l}(\lambda_{l-1},d\lambda_{l})f(\lambda_{(k+i-p+1)^{+}:k+i})\leq1.\]
 and conclude by noting that
 \begin{eqnarray*}
\tilde{\gamma}_{k+1:k+i}^{p}\tilde{R}_{k}^{p}\mu(1) & = & \int(\tilde{R}_{k}^{p}\mu)(d\lambda_{(k-p+1)^{+}:k})\prod_{l=k+1}^{k+i}Q_{l}(\lambda_{l-1},d\lambda_{l})\tilde{\Psi}_{l}^{p}(\lambda_{(l-p+1)^{+}:l})\\
 & \geq & \frac{1}{a_{k+1}\ldots a_{k+i}}\end{eqnarray*}
 since $\tilde{R}_{k}^{p}\mu$ is a probability measure.
\end{proof}

\section{Mixing and contraction properties of the truncated filter\label{sec:mix-truncated}}

The truncated filter may be interpreted as a standard filter based
on Markov chain $\tilde{\Lambda}_{n}^{p}=\Lambda_{(n-p+1)^{+}:n}$. This
insight allows us to establish the contraction properties of the truncated
filter.
\begin{lem}
\label{lem:melange} One has:
 \[
h(\tilde{R}_{k+1:k+p}^{p}\mu,\tilde{R}_{k+1:k+p}^{p}\mu')\leq\frac{1}{\tilde{\varepsilon}_{k+1,p}^{2}}\Vert\mu-\mu'\Vert_{TV}\]
 and
 \[
h(\tilde{R}_{k+1:k+p}^{p}\mu,\tilde{R}_{k+1:k+p}^{p}\mu')\leq\tilde{\rho}_{k+1,p}h(\mu,\mu')\]
 where
\[
\tilde{\varepsilon}_{k,p}^{2}=\frac{\varepsilon_{k}^{2}}{(a_{k}\ldots a_{k+p-2})(b_{k}\ldots b_{k+p-2})},\qquad\tilde{\rho}_{k,p}=\frac{1-\tilde{\varepsilon}_{k,p}^{2}}{1+\tilde{\varepsilon}_{k,p}^{2}},\]
for all $k\geq0$, and all $\mu$ ,$\mu'\in\mathcal{P}(E^{(k+1)\wedge p})$.
\end{lem}
Note $\tilde{\varepsilon}_{k,n}$ must be interpreted as a mixing
coefficient, and $\tilde{\rho}_{k,p}$ as a Birkhoff contraction coefficient.
\begin{proof}
Using Hypothesis 3, one has:
 \begin{eqnarray*}
\lefteqn{Q_{k+p}\tilde{\gamma}_{k+1:k+p-1}\mu}\\
 & = & \int\mu(d\lambda_{(k-p+1)^{+}:k})\prod_{i=k+1}^{k+p}Q_{i}(\lambda_{i-1},d\lambda_{i})\prod_{i=k+1}^{k+p-1}\left[\tilde{\Psi}_{i}^{p}(\lambda_{(i-p+1)^{+}:i})\right]\\
 & \leq & b_{k+1}\dots b_{k+p-1}\int\mu(d\lambda_{(k-p+1)^{+}:k})\prod_{i=k+1}^{k+p}Q_{i}(\lambda_{i-1},d\lambda_{i})\\
 & \leq & \frac{b_{k+1}\dots b_{k+p-1}}{\varepsilon_{k+1}}\tilde{\xi}_{p}(d\lambda_{k+1:k+p})\end{eqnarray*}
 where $\tilde{\xi}_{p}$ stands for the following reference measure:
\[
\tilde{\xi}_{p}(d\lambda_{k+1:k+p})=\xi(d\lambda_{k+1})\prod_{i=k+2}^{k+p}Q_{i}(\lambda_{i-1},d\lambda_{i}).\]
 One shows similarly that
 \[
Q_{k+p}\tilde{\gamma}_{k+1:k+p-1}\mu\geq\frac{\varepsilon_{k+1}}{a_{k+1}\dots a_{k+p-1}}\tilde{\xi}_{p}(d\lambda_{k+1:k+p}).\]
 Hence kernel $Q_{k+p}\tilde{\gamma}_{k+1:k+p-1}\mu$ is mixing, with
mixing coefficient $\tilde{\varepsilon}_{k+1,p}$.

Following Lemma 3.4 in \citet{legland2004sau},
 \begin{eqnarray*}
h(\tilde{R}_{k+1:k+p}^{p}\mu,\tilde{R}_{k+1:k+p}^{p}\mu') & = & h(Q_{k+p}\tilde{\gamma}_{k+1:k+p-1}\mu,Q_{k+p}\tilde{\gamma}_{k+1:k+p-1}\mu')\\
 & \leq & \frac{1}{\tilde{\varepsilon}_{k+1,p}^{2}}\tv{\mu-\mu'}\end{eqnarray*}
 using the scale invariance property of the Hilbert metric. Similarly,
according to Lemma 3.9 in the same paper:
 \begin{eqnarray*}
h(\tilde{R}_{k+1:k+p}^{p}\mu,\tilde{R}_{k+1:p}^{p}\mu') & = & h(Q_{k+p}\tilde{\gamma}_{k+1:k+p-1}\mu,Q_{k+p}\tilde{\gamma}_{k+1:k+p-1}\mu')\\
 & \leq & \left(\frac{1-\tilde{\varepsilon}_{k+1,p}^{2}}{1+\tilde{\varepsilon}_{k+1,p}^{2}}\right)h(\mu,\mu').
\end{eqnarray*}

\end{proof}

\section{\label{sec:Propagation}Propagation of  truncation error}

We establish first the two following lemmas.
\begin{lem}
\label{lem:tele1} Let $\bar{R}_{n}:E^{n\wedge p}\rightarrow\mathcal{P}(E^{(n+1)\wedge p})$
be a sequence of (possibly random) probability kernels such that for
all $n\geq1$ and $\mu\in\mathcal{P}(E^{n\wedge p})$,

\[
\sup_{f:\Vert f\Vert_{\infty}=1}\mathbb{E}\left\{ \left|\langle\tilde{R}_{n}^{p}\mu-\bar{R}_{n}\mu,f\rangle\right|\right\} \leq\delta_{n}\ ,\]
where the expectation is w.r.t. the randomness of $\bar{R}_n$, 
then, for all $n\geq1$ and all $\zeta\in\mathcal{P}(E)$, one has
\[
\sup_{f:\Vert f\Vert_{\infty}=1}\mathbb{E}\left\{ \left|\langle\tilde{R}_{1:n}^{p}\zeta-\bar{R}_{1:n}\zeta,f\rangle\right|\right\} \leq\frac{8}{\log(3)}\sum_{i=1}^{n}\left(\frac{\delta_{i}}{\tilde{\varepsilon}_{i+1}^{2}\tilde{\varepsilon}_{i+p+1}^{2}}\prod_{j=2}^{\lfloor\frac{n-i}{p}\rfloor-1}\tilde{\rho}_{i+jp+1,p}\right)\]
where $\bar{R}_{1:n}\zeta=\bar{R}_{n}\ldots\bar{R}_{1}\zeta$, and
with the convention that empty products equal one. \end{lem}
\begin{proof}
The following difference can be decomposed into a telescopic sum:
\[
\tilde{R}_{1:n}^{p}\zeta-\bar{R}_{1:n}\zeta=\sum_{i=1}^{n}\left(\tilde{R}_{i+1:n}^{p}\tilde{R}_{i}^{p}\bar{R}_{1:i-1}\zeta-\tilde{R}_{i+1:n}^{p}\bar{R}_{i}\bar{R}_{1:i-1}\zeta\right).\]
We fix the integers $i$, $n$, and consider some arbitrary test function
$f$. For $i\geq n-2p$, one may apply Lemma \ref{lem:32}:
 \begin{eqnarray*}
\lefteqn{\sup_{f:\Vert f\Vert_{\infty}=1}\mathbb{E}
\left\{
  \left|\langle\tilde{R}_{i+1:n}^{p}\tilde{R}_{i}^{p}\bar{R}_{1:i-1}\zeta
-\tilde{R}_{i+1:n}^{p}\bar{R}_{i}\bar{R}_{1:i-1}\zeta,f\rangle\right|\right\} }\\
 & \leq & 2(a_{i+1}\dots a_{n})(b_{i+1}\dots b_{n})\delta_{i}\\
 &  \leq & \frac{8}{\log(3)}\frac{\delta_{i}}
{\tilde{\varepsilon}_{i+1,p}^{2}\tilde{\varepsilon}_{i+p+1,p}^{2}}
\end{eqnarray*}
 since $\varepsilon_{n}\leq1$, $a_{n}\geq1$ and $b_{n}\geq1$ for
all $n$.

For $i<n-2p$, let $k=\lfloor(n-i)/p\rfloor$, then, using Lemma \ref{lem:melange},
Equations (\ref{Eq:borne1}) to \eqref{Eq:inv1}
one has \begin{eqnarray*}
\lefteqn{\left|\langle\tilde{R}_{i+1:n}^{p}\tilde{R}_{i}^{p}\bar{R}_{1:i-1}\zeta-\tilde{R}_{i+1:n}^{p}\bar{R}_{i}\bar{R}_{1:i-1}\zeta,f\rangle\right|}\\
 & \leq & \left\Vert \tilde{R}_{i+1:n}^{p}\tilde{R}_{i}^{p}\bar{R}_{1:i-1}\zeta-\tilde{R}_{i+1:n}^{p}\bar{R}_{i}\bar{R}_{1:i-1}\zeta\right\Vert _{TV}\\
 & \leq & \frac{2}{\log(3)}h\left(\tilde{R}_{i+1:i+kp}^{p}\tilde{R}_{i}^{p}\bar{R}_{1:i-1}\zeta,\tilde{R}_{i+1:i+kp}^{p}\bar{R}_{i}\bar{R}_{1:i-1}\zeta\right)\\
 & \leq & \frac{2}{\log(3)\tilde{\varepsilon}_{i+p+1,p}^{2}}\times\prod_{j=2}^{k-1}\tilde{\rho}_{i+jp+1,p}\times\left\Vert \tilde{R}_{i+1:i+p}^{p}\nu-\tilde{R}_{i+1:i+p}^{p}\nu'\right\Vert _{TV}\end{eqnarray*}
 where $\nu=\tilde{R}_{i}^{p}\bar{R}_{1:i-1}\zeta$, $\nu'=\bar{R}_{i}\bar{R}_{1:i-1}\zeta$.
Applying (7) p. 160 of \citet{legland2004sau}, one gets
 \begin{eqnarray*}
\left\Vert \tilde{R}_{i+1:i+p}^{p}\nu-\tilde{R}_{i+1:i+p}^{p}\nu'\right\Vert _{TV} & \leq & 2\,\frac{\Vert\tilde{\gamma}_{i+1:i+p}^p\nu-\tilde{\gamma}_{i+1:i+p}^p\nu'\Vert_{TV}}{\tilde{\gamma}^p_{i+1:i+p}\nu(1)}.\end{eqnarray*}
 where, using the same calculations as in Lemma \ref{lem:melange},
\[
\tilde{\gamma}^p_{i+1:i+p}\nu(1)\geq\frac{\varepsilon_{i+1}}{a_{i+1}\dots a_{i+p}}\]
 and
 \begin{eqnarray*}
\lefteqn{\mathbb{E}\left[\left\Vert \tilde{\gamma}^p_{i+1:i+p}\nu-\tilde{\gamma}^p_{i+1:i+p}\nu'\right\Vert _{TV}\right]}\\
 & = & \mathbb{E}\left[\int_{x'\in E^{p}}\left|\int_{x\in E^{p}}(\nu-\nu')(dx)\tilde{\gamma}^p_{i+1:i+p}(x,dx')\right|\right]\\
 & \leq & \left[\sup_{x\in E^{p}}\int_{x'\in E^{p}}\tilde{\gamma}^p_{i+1:i+p}(x,dx')\right]\left[\sup_{\phi:\Vert\phi\Vert_{\infty}=1}\mathbb{E}(|\langle\nu-\nu',\phi\rangle|)\right]\\
 & \leq & \frac{b_{i+1}\ldots b_{i+p}}{\varepsilon_{i+1}}\left[\sup_{\phi:\Vert\phi\Vert_{\infty}=1}\mathbb{E}(|\langle\nu-\nu',\phi\rangle|)\right]\end{eqnarray*}
 which ends the proof.\end{proof}
\begin{lem}
\label{lem:tele2} For all $n\geq1$ and all $\zeta\in\mathcal{P}(E)$,
one has
$$\left\Vert \tilde{R}_{1:n}^{p}\zeta-H_{n}^{p}R_{1:n}\zeta\right\Vert _{TV}\leq\frac{4\tau^{p}}{\log3}
\left\{ \sum_{i=1}^{n}\frac{\phi_{i}}{\tilde{\varepsilon}_{i+1,p}^{2}}\prod_{j=1}^{\left\lfloor (n-i)/p\right\rfloor -1}\tilde{\rho}_{i+jp+1,p}\right\}
$$
with the convention that empty sums equal zero, and empty products
equal one.\end{lem}
\begin{proof}
One has:
\[
\tilde{R}_{1:n}^{p}\zeta-H_{n}^p
R_{1:n}\zeta=\sum_{i=1}^{n}\left(\tilde{R}_{i+1:n}^{p}\tilde{R}_{i}^{p}H_{i-1}^{p
 } R_{1:i-1}\zeta-\tilde{R}_{i+1:n}^{p}H_{i}^p R_{1:i}\zeta\right)\]
 For $i\leq n-p$, let $k=\lfloor(n-i)/p\rfloor$, then according
to Lemma \ref{lem:melange}: \begin{eqnarray*}
 &  & \left\Vert \tilde{R}_{i+1:n}^{p}\tilde{R}_{i}^{p}H_{i-1}^{p}R_{1:i-1}\zeta-\tilde{R}_{i+1:n}^{p}\tilde{R}_{i+1}^{p}H_{i}^{p}R_{1:i}\zeta\right\Vert _{TV}\\
 &  & \quad\leq\frac{2}{\log3}h\left(\tilde{R}_{i+1:i+kp}^{p}\tilde{R}_{i}^{p}H_{i-1}^{p}R_{1:i}\zeta,\tilde{R}_{i+1:i+kp}^{p}H_{i}^{p}R_{1:i}\zeta\right)\\
 &  & \quad\leq\frac{2}{\log(3)\tilde{\varepsilon}_{i+1,p}^{2}}\prod_{j=1}^{k-1}\tilde{\rho}_{i+jp+1,p}\left\Vert \tilde{R}_{i}^{p}H_{i-1}^{p}R_{1:i-1}\zeta-H_{i}^{p}R_{1:i}\zeta\right\Vert _{TV}\end{eqnarray*}
 and ones concludes using Lemma \ref{lem:loc2}. For $i>n-p$, one
can apply Lemma \ref{lem:loc2} directly.
\end{proof}

\section{\label{sec:Coupling}Coupling of particle approximations}

We now introduce two interactive particle systems: the first particle
system approximates the true filter, and is equivalent to the type
of particle algorithms studied in this paper, and the second particle
system approximates the truncated filter, and corresponds to an artificial
algorithm that would not be implemented in practice. We work out a
way of coupling both particle systems in order to evaluate the distance
between the two (in a sense that is made clear below).

We define, for $n\geq0$,
 \begin{eqnarray*}
\bar{Q}_{n,p}\left(\lambda_{(n-p)^{+}:n-1},d\lambda'_{(n-p+1)^{+}:n}\right) 
& = & \delta_{\lambda_{(n-p+1)^{+}:n-1}}(d\lambda'_{(n-p+1)^{+}:n-1})\\
 &  & ~~~\times Q_{n}(\lambda_{n-1},d\lambda'_{n}),
\end{eqnarray*}
 \[
\bar{Q}_{n}(\lambda_{0:n-1},d\lambda'_{0:n})=\delta_{\lambda_{0:n-1}}(d\lambda'_{0:n-1})\times Q_{n}(\lambda_{n-1},d\lambda'_{n})\ .\]
 We define $\forall\nu\in\mathcal{M}_{+}(E^{n+1})$, $\forall$ measurable
$f:E^{n+1}\rightarrow\mathbb{R}^{+}$, $\forall\nu'\in\mathcal{M}_{+}(E^{p})$,
$\forall$ measurable $g:E^{(n+1)\wedge p}\rightarrow\mathbb{R}^{+}$,
\[
\Psi_{n}.\nu(f)=\frac{\nu(\Psi_{n}f)}{\nu(\Psi_{n})}\ ,
\tilde{\Psi}_{n}^{p}.\nu'(g)=\frac{\nu'(\tilde{\Psi}_{n}^{p}g)}{\nu'(\Psi_{n}^{p})}\ .\]
 For any measurable space $(E',\Omega')$ and any measure $\mu'\in\mathcal{P}(E')$,
we can take $Z_{1},Z_{2},\dots$ i.i.d. of law $\mu'$ and define
the random empirical measure, for $N\geq 1$,
 \[
S^{N}(\mu')=\frac{1}{N}\sum_{i=1}^{N}\delta_{Z_{i}}\ .\]
 Notice that, as the $Z_{1},Z_{2},\dots$ are only given in law, we
only define $S^{N}(\mu)$ in law. We define the random operators $R_{n}^{N}$,
$\tilde{R}_{n}^{p,N}$ ($\forall n$) by: $\forall\mu\in\mathcal{P}(E^{n})$,
$R_{n}^{N}\mu$ is a random weighted empirical measure such that
\[
R_{n}^{N}\mu=\Psi_{n}.S^{N}(\bar{Q}_{n}\mu)\ .\]
 Similarly, $\forall\mu'\in\mathcal{P}(E^{p\wedge n})$, $\tilde{R}_{n}^{p,N}\mu'$
is a random weighted empirical measure such that
\begin{equation}
\label{eq:inlaw}
\tilde{R}_{n}^{p,N}\mu'=\tilde{\Psi}_{n}^{p}.S^{N}(\bar{Q}_{n,p}\mu')\ .
\end{equation}
 As pointed above, $R_{n}^{N}\mu$ and $\tilde{R}_{n}^{p,N}\mu'$
are only defined in law. Since $\zeta$ denotes the probability density
of the first state $\Lambda_{0}$, the particle system with $N$ particles
approximating the true filter at time $n$ is defined by
 \[
 R_{n}^{N}R_{n-1}^{N}\dots R_{1}^{N}\zeta,\] 
and the particle system
 with $N$ particles approximating the truncated filter at time $n$ is
 defined by
\[
\tilde{R}_{n}^{p,N}\tilde{R}_{n-1}^{p,N}\dots\tilde{R}_{1}^{p,N}\zeta.\]

\begin{lem}
\label{lem:coupling}There exists a coupling such that, for all $k\geq1$
and  $\mu\in\mathcal{P}(E^k)$:
 \[
\sup_{f:\Vert f\Vert_{\infty}\leq1}\mathbb{E}\left(\left|\langle\tilde{R}_{k}^{p,N}H_{k-1}^{p}\mu-H_{k}^{p}R_{k}^{N}\mu,f\rangle\right|\right)\leq\phi_{k}\tau^{p}.\]

As $H_{k}^{p}R_{k}^{N}\mu$ and $\tilde{R}_{k}^{p,N}H_{k}^{p}\mu$
are defined to be random variables with such and such law, the term
``coupling'' means that we can define a
random variable $(H_{k}R_{k}^{N}\mu,\tilde{R}_{k}^{p,N}H_{k}\mu)$
with the desired marginals. \end{lem}
\begin{proof}
To prove the above result, we produce a coupling between the two random
measures $\tilde{R}_{k}^{p,N}H_{k-1}^{p}\mu$ and $H_{k}^{p}R_{k}^{N}\mu$.
Let \[
\bar{\Psi}_{n}(\lambda_{0:n})=\tilde{\Psi}_{n}^{p}(\lambda_{(n-p+1)^{+}:n}),\]
so that, for $\mu\in\mathcal{P}(E^{k})$, and using \eqref{eq:inlaw},
one has
 \[
\tilde{R}_{k}^{p,N}H_{k-1}^{p}\mu=H_{k}^{p}(\bar{\Psi}_{k}.(S^{N}(\bar{Q}_{k}\mu))\,\]
in the sense that both sides define the same distribution. 
Let $\chi_{1},\dots,\chi_{N}$ i.i.d. $\sim\mu\bar{Q}_{k}$, where $\chi_{i}$
is a vector $\lambda_{0:k,i}$, for $i=1,\ldots,N$, and $\tilde{\chi}_{i}$
denotes its projection on the $p$ last components, $\tilde{\chi}_{i}=\lambda_{(k-p+1)^{+}:k,i},$
then \[
\frac{1}{\sum_{j=1}^{N}\Psi_{k}(\chi_{j})}\sum_{i=1}^{N}\Psi_{k}(\chi_{i})\delta_{\hat{\chi_{i}}}\text{ has same law as }H_{k}^{p}R_{k}^{N}\mu\]
 and \[
\frac{1}{\sum_{j=1}^{N}\bar{\Psi}_{k}(\chi_{j})}\sum_{i=1}^{N}\bar{\Psi}_{k}(\chi_{i})\delta_{\tilde{\chi}_{i}}\text{ has same law as }\tilde{R}_{k}^{p,N}H_{k-1}^{p}\mu\ .\]
 For any $f$ such that $\Vert f\Vert_{\infty}\leq1$ (using a classical
result on empirical measures):
 \begin{eqnarray*}
 &  & |\langle\tilde{R}_{k}^{p,N}H_{k-1}^{p}\mu-H_{k}^{p}R_{k}^{N}\mu,f\rangle|\\
 &  & \qquad\leq\frac{1}{2}\sum_{i=1}^{N}\left|\frac{\Psi_{k}(\chi_{i})}{\sum_{j=1}^{N}\Psi_{k}(\chi_{j})}-\frac{\bar{\Psi}_{k}(\chi_{i})}{\sum_{j=1}^{N}\bar{\Psi}_{k}(\chi_{j})}\right|\\
 &  &
 \qquad\leq\frac{1}{2}\sum_{i=1}^{N}\left(\left|\frac{\Psi_{k}(\chi_{i})-\bar{\Psi}_{k}(\chi_{i})}{\sum_{j=1}^{N}\Psi_{k}(\chi_{j})}\right|
+\left|\frac{\bar{\Psi}_{k}(\chi_{i})\sum_{j=1}^{N}(\bar{\Psi}_{k}(\chi_{j})-\Psi_{k}(\chi_{j}))}
  {\left\{\sum_{j=1}^{N}\Psi_{k}(\chi_{j})\right\}
    \left\{\sum_{j=1}^{N}\bar{\Psi}_{k}(\chi_{j})\right\}}\right|\right)\\
 &  &
 \qquad\leq\frac{\phi_{k}\tau^{p}}{2}\sum_{i=1}^{N}\left(\frac{\Psi_{k}(\chi_{i})}{\sum_{j=1}^{N}\Psi_{k}(\chi_{i})}+\frac{\bar{\Psi}_{k}(\chi_{i})\sum_{j=1}^{N}\bar{\Psi}_{k}(\chi_{i})\wedge{\Psi}(\chi_{i})}{\left\{
       \sum_{j=1}^{N}\Psi(\chi_{j})\right\} \left\{
       \sum_{j=1}^{N}\bar{\Psi}(\chi_j)\right\} }\right)\\
 &  & \qquad\leq\phi_{k}\tau^{p}\,\end{eqnarray*}
using Hypothesis 2, from which we deduce the result.
\end{proof}

\section{\label{sec:Main}Main result}

We are now able to derive estimates of the error
\[
\mathcal{E}_{n,N}^{p}(y_{1:n})=\sup_{f:\Vert f\Vert_{\infty}=1}\mathbb{E}_{N}\left(\left|\langle H_{n}^{p}R_{1:n}\zeta-H_{n}^{p}R_{1:n}^{N}\zeta,f\rangle\right|\ |Y_{0:n}=y_{0:n}\right)\]
 induced by the particle approximation of the true filter, for the
marginal filtering distribution of the $p$ last states, provided
$p\leq n$. The expectation $\mathbb{E}_{N}$ is with respect to the
randomness of the $N$ particles, and the functions $f$ are $E^{(n+1)\wedge p}\rightarrow\mathbb{R}$.
Note that $\mathcal{E}_{n,N}^{p}(y_{1:n})$ is by construction an
increasing function of $p$. 
\begin{thm}
\label{prop:maj} For any $\zeta\in\mathcal{P}(E)$, and any test function
w.r.t. $E^{(n+1)\wedge p}$,

\begin{equation}
\mathcal{E}_{n,N}^{p}(y_{1:n})\leq\frac{4}{\log3}\sum_{i=1}^{n}\frac{\delta_{i}}{\tilde{\varepsilon}_{i+1,p}^{2}\tilde{\varepsilon}_{i+p+1,p}^{2}}\prod_{j=2}^{\lfloor(n-i)/p\rfloor-1}\tilde{\rho}_{i+jp+1}\label{eq:error_bound}\end{equation}
where \[
\delta_{i}=3\tau^{p}\phi_{i}+\frac{4a_{i}b_{i}}{\sqrt{N}}.\]
\end{thm}
\begin{proof}
We first study the following local error, for $\mu\in\mathcal{P}(E^{n})$,
\[
\sup_{f:\left\Vert f\right\Vert _{\infty}=1}\mathbb{E}_N
\left[\left|\left\langle
      \tilde{R}_{n}^{p}H_{n-1}^{p}\mu-H_{n}^{p}R_{n}^{N}\mu,f\right\rangle
  \right|\ \bigm\vert Y_{0:n}=y_{0:n}\right]\]
 where the difference of operators can de decomposed into:
 \[
\tilde{R}_{n}^{p}H_{n-1}^{p}-H_{n}^{p}R_{n}^{N}=\left(\tilde{R}_{n}^{p}H_{n-1}^{p}-\tilde{R}_{n}^{p,N}H_{n-1}^{p}\right)+\left(\tilde{R}_{n}^{p,N}H_{n-1}^{p}-H_{n}^{p}R_{n}^{N}\right).\]
 To bound the first term, one may use (25) p. 162 of \citet{legland2004sau},
for $\nu=H_{n-1}^{p}\mu$ and Hypothesis 3:
 \[
 \mathbb{E}_{N}\left[\left|\left\langle
       \tilde{R}_{n}^{p}\nu-\tilde{R}_{n}^{p,N}H_{n-1}^{p}\nu,f\right\rangle
   \right|\right]\leq\frac{2a_{n}b_{n}}{\sqrt{N}}\] 
and, for the
 second term, one may apply Lemma \ref{lem:coupling}:
\[
 \mathbb{E}_N\left[\left|\left\langle
       \tilde{R}_{n}^{p,N}H_{n-1}^{p}\mu-H_{n}^{p}R_{n}^{N}\mu,f\right\rangle
   \right|\right]\leq\phi_{n}\tau^{p}\] 
so that
\[ \sup_{f:\left\Vert
     f\right\Vert _{\infty}=1}\mathbb{E}_{N}\left[\left|\left\langle
       \tilde{R}_{n}^{p}\mu-H_{n}^p R_{n}^{N}\mu,f\right\rangle
   \right|\right]\leq\delta_{n}'\] 
for $\delta_{n}'=2a_{n}b_{n}/\sqrt{N}+\phi_{n}\tau^{p}$. This local error
 is propagated using Lemma \ref{lem:tele1}:
 \begin{eqnarray*}
   \lefteqn{\mathbb{E}_{N}\left[\left|\left\langle \tilde{R}_{1:n}^{p}\zeta-H_{n}^{p}R_{1:n}^{N}\zeta,f\right\rangle \right|\right]}\\
   & \leq &
   \frac{8}{\log(3)}\sum_{i=1}^{n}\left(\frac{\delta_{i}'}{\tilde{\varepsilon}_{i+1}^{2}\tilde{\varepsilon}_{i+p+1}^{2}}
\prod_{j=2}^{\lfloor\frac{n-i}{p}\rfloor-1}\tilde{\rho}_{i+jp+1,p}\right).
\end{eqnarray*}

To conclude, one may decompose the global error as follows:
 \[
H_{n}^{p}R_{1:n}\zeta-H_{n}^{p}R_{1:n}^{N}\zeta=\left(H_{n}^{p}R_{1:n}\zeta-\tilde{R}_{1:n}^{p}\zeta\right)+\left(\tilde{R}_{1:n}^{p}\zeta-H_{n}^{p}R_{1:n}^{N}\zeta\right).\]
 where the second term is bounded above, and the first term is directly
bounded using Lemma \ref{lem:tele2}.
\end{proof}
Since $p$ is an arbitrary parameter, one may minimise the error bound
with respect to $p$. For instance, one has the following result for
time-uniform estimates. As noted above, the error $\mathcal{E}_{n,N}^{p}(y_{1:n})$
is an increasing function of $p$, so the bound below applies a fortiori
to $\mathcal{E}_{n,N}^{1}(y_{1:n})$, the particle error corresponding
to the marginal filtering distribution of the last state $\Lambda_{n}$. 
\begin{cor}
\label{cor:time-unif}If there exists constants $c$, $\varepsilon$,
$\phi>0$ such that, almost surely, $a_{n}b_{n}\leq c$, $\varepsilon_{n}\ge\varepsilon$,
and $\phi_{n}\leq\phi$, then, provided $\tau c^{3}<1$, the particle
error is bounded almost surely as follows:
 \[
\mathcal{E}_{n,N}^{p}(y_{1:n})\leq C\left\{ \log(N)+D\right\} \left(\frac{1}{\sqrt{N}}\right)^{1+3\log c/\log\tau},\]
for $N$ large enough, where
 \[
C=\frac{16}{\varepsilon^{6}c^2}
\left(\frac{-1}{\log\tau}\right)
\left(\frac{4c}{3\phi}\right)^{3\log  c/\log\tau},
\quad D=2\log(3\phi/4c\tau),\]
and
\begin{equation}
p=\left\lceil \log\left\{ \frac{4c}{3\phi\sqrt{N}}\right\} \slash\log\tau\right\rceil .\label{eq:pfN}\end{equation}
\end{cor}
\begin{proof}
Under these conditions, the RHS of (\ref{eq:error_bound}) is smaller
than or equal to:
\begin{eqnarray}
\mathcal{E}_{n,N}^{p}(y_{1:n}) & \leq & \frac{4}{\log3}\frac{c^{2(p-2)}}{\varepsilon^{4}}\left(3\phi\tau^{p}+\frac{4c}{\sqrt{N}}\right)\sum_{i=1}^{n}\left(1-\varepsilon^{2}c^{-(p-2)}\right)^{\lfloor(n-i)/p\rfloor-2}\label{eq:error_cor}\\
 & \leq & \frac{4}{\log3}\frac{c^{2(p-2)}}{\varepsilon^{4}}\left(3\phi\tau^{p}+\frac{4c}{\sqrt{N}}\right)\sum_{i=0}^{n-1}\left(1-\varepsilon^{2}c^{-(p-2)}\right)^{i/p-1}\nonumber \\
 & \leq & \frac{4}{\log3}\frac{c^{2(p-2)}}{\varepsilon^{4}}\left(3\phi\tau^{p}+\frac{4c}{\sqrt{N}}\right)\frac{\left(1-\varepsilon^{2}c^{-(p-2)}\right)^{-1}}{1-\left(1-\varepsilon^{2}c^{-(p-2)}\right)^{1/p}}\nonumber \\
 & \leq & \frac{4c^{3(p-2)}}{\varepsilon^{6}}\left(3\phi\tau^{p}+\frac{4c}{\sqrt{N}}\right)p\nonumber 
\end{eqnarray}
for $p$ large enough, since $(1-x)^{a}\leq1-ax$ for $a\in(0,1)$,
$x\in(0,1)$, so, provided $c^{3}\tau<1$, one may take $p$ as in
\eqref{eq:pfN}, 
which gives:
 \[
\mathcal{E}_{n,N}^{p}(y_{1:n})\leq
\frac{32}{\varepsilon^{6}c^2}
\left(\frac{4c}{3\phi}\right)^{3 \frac{\log c}{\log \tau}}
\left(\frac{\log
    N+2\log(3\phi/4c)}{-2\log\tau}+1\right)
\left(\frac{1}{\sqrt{N}}\right)^{1+\frac{3\log c}{\log\tau}}\]
and conclude. \[
\]
\end{proof}
Obviously, this is a qualitative result, in that there are many practical
models where such time-uniform, deterministic bounds are not available.
For specific models, one may be able instead to use (\ref{eq:error_bound})
in order to establish the asymptotic stability of the expected particle
error, where the expectation is with respect to observed process $(Y_{n})$.
We provide an example of this approach in Section \ref{sec:applis}.

\section{Applications to practical models\label{sec:applis}}

In this section, we apply our general result to three practical models.
We keep the same settings and notations, i.e. the observed process
$(Y_{n})$ admits some probability distribution conditional on the
path $\Lambda_{0:n}=\lambda_{0:n}$ of a Markov chain $(\Lambda_{n})$,
with initial distribution $\zeta$ and Markov transition $Q_{n}$, which
fulfil Hypothesis 1, see Section \ref{sec:model}. We derive conditions on
the model parameters that ensure asymptotic stability of the particle
error; in particular, these conditions imply that Hypotheses 2 and
3 are verified.

We state the following trivial result for further reference. Let $(f,g)$
a pair of probability densities $(f,g)$ on $E$, then:
 \begin{eqnarray}
\lefteqn{\forall x\in E,\,|\log\left\{ f(x)\right\} -\log\left\{ g(x)\right\} |\leq c}\nonumber \\
 & \Rightarrow & \forall x\in E,\,|f(x)-g(x)|\leq(e^{c}-1)\left\{ f(x)\wedge g(x)\right\} \label{ineq:logdens}
\end{eqnarray}
for $c\geq0$.

\subsection{GARCH Mixture model}

We assume that the observed process is such that
 \[
 Y_{n}=\sigma_{n}(\Lambda_{0:n})Z_{n},\quad n\geq1,\] 
where the $Z_{n}$'s are i.i.d. $\mathcal{N}(0,1)$ random variables, and the
 variance function $\sigma_{n}^{2}$ is defined recursively, for
 $n\geq1$:
 \begin{equation}
   \sigma_{n}^{2}(\lambda_{0:n})=\alpha(\lambda_{n})+\beta(\lambda_{n})Y_{n}^{2}+\gamma(\lambda_{n})\sigma_{n-1}^{2}(\lambda_{0:n-1})\label{eq:sigman}
\end{equation}
 and $\sigma_{0}^{2}(\lambda_{0})=\alpha(\lambda_{0})/\left\{
   1-\gamma(\lambda_{0})\right\} ,$ where $\alpha$, $\beta$ and
 $\gamma$ are $E\rightarrow\mathbb{R}^{+}$ functions. Conditional on
 $\Lambda_{0:n}$, $(Y_n)$ is a GARCH (generalised
 autoregressive conditional heteroskedasticity) process 
 \citep{bollerslev1986gac}; see \citet{Chopin:change}
for a finance application of such a GARCH mixture model. 

The potential functions equal
 \[
\Psi_{n}(\lambda_{0:n})=\frac{1}{\sqrt{2\pi\sigma_{n}^{2}(\lambda_{0:n}})}\exp\left\{ -\frac{y_{n}^{2}}{2\sigma_{n}^{2}(\lambda_{0:n})}\right\},
\]
for $\lambda_{0:n}\in E^{n+1}$, and $(\Lambda_{n})$ is a Markov
process, with Markov kernels $Q_{n}$, which satisfy Hypothesis 1. 

The functions $\alpha$, $\beta$ and $\gamma$ are assumed to be
bounded as follows:\[
0<\alpha_{\min}\leq\alpha(\lambda)\leq\alpha_{\max},\quad0\le\beta_{\min}\leq\beta(\lambda)\leq\beta_{\max}<1,\]
\[
0\leq\gamma_{\min}\leq\gamma(\lambda)\leq\gamma_{\max}<1.\]

We first consider the case where $\beta(\lambda)=0$ for all $\lambda\in E$.
As mentioned in the introduction, this simplified model can be interpreted
as a standard hidden Markov model, with observed process $(Y_{n})$,
and Markov chain $\left(\kappa_{n}\right)=\left(\Lambda_{n},\sigma_{n}^{2}(\Lambda_{0:n})\right)$.
However, since $\sigma_{n}^{2}(\Lambda_{0:n})$ is a deterministic
function of $\sigma_{n-1}^{2}(\Lambda_{0:n-1})$ and $\lambda_{n}$,
it does not have mixing or similar properties that are usually required
to obtain estimates of the particle error. Instead, analysing this
model as a Feynman-Kac flow with iterative, path-dependent potential
functions make it possible to derive such estimates. 
\begin{lem}
For the simplified model described above (with $\beta=0$), the expected
particle error of the corresponding particle approximation is uniformly
stable in time, i.e. there exists constants $C$, $D$, such that
\[
\mathbb{E}\left[\mathcal{E}_{n,N}^{p}(Y_{1:n})\right]\leq C\left\{ \log(N)+D\right\} \left(\frac{1}{\sqrt{N}}\right)^{1+3\log c/\log\tau},\]
 where $p$ is given by (\ref{eq:pfN}), provided $\iota<2$ and $\tau c^{3}<1$,
where $\tau=\gamma_{\max}$, $c=\left(2/\iota-1\right)^{-1/2}$, and
\[
\iota=\frac{\alpha_{\max}\left(1-\gamma_{\min}\right)}{\alpha_{\min}\left(1-\gamma_{\max}\right)}.\]
\end{lem}
\begin{proof}
From (\ref{eq:sigman}), one sees the process $\sigma_{n}^{2}$ is
bounded, $\sigma_{\min}^{2}\leq\sigma_{n}^{2}(\lambda_{0:n})\leq\sigma_{\max}^{2}$
for all $\lambda_{0:n}\in E^{n+1}$, where
 \[
\sigma_{\min}^{2}=\frac{\alpha_{\min}}{1-\gamma_{\min}},\quad\sigma_{\max}^{2}=\frac{\alpha_{\max}}{1-\gamma_{\max}}.\]
 so, for a given sequence observations $y_{1:n}$, Hypothesis 3 is
verified with:
\[
\frac{1}{a_{n}}=\frac{1}{\sqrt{2\pi\sigma_{\max}^{2}}}\exp\left\{ -\frac{y_{n}^{2}}{2\sigma_{\min}^{2}}\right\} ,\quad b_{n}=\frac{1}{\sqrt{2\pi\sigma_{\min}^{2}}}\exp\left\{ -\frac{y_{n}^{2}}{2\sigma_{\max}^{2}}\right\} ,\]
provided the truncated potential is taken as:
\[
\tilde{\Psi}_{n}^{p}(\lambda_{n-p+1:n})=\Psi_{n}(z,\ldots,z,\lambda_{n-p+1:n})\]
where $z$ is an arbitrary element of $E$. For Hypothesis 2, one
has, for any $\lambda_{0:n},\lambda_{0:n}'\in E^{(n+1)}$ such that
$\lambda_{(n-p+1)^{+}:n}=\lambda_{(n-p+1)^{+}:n}'$ :
 \begin{eqnarray*}
\left|\log\Psi_{n}(\lambda_{0:n})-\log\Psi_{n}(\lambda_{0:n}')\right| & \leq & \frac{1}{2}\left|\log\sigma_{n}^{2}(\lambda_{0:n})-\log\sigma_{n}^{2}(\lambda_{0:n}')\right|\\
 &  & +\frac{y_{n}^{2}}{2}\left|\frac{1}{\sigma_{n}^{2}(\lambda_{0:n})}-\frac{1}{\sigma_{n}^{2}(\lambda_{0:n}')}\right|\\
 & \leq & \frac{\sigma_{\min}^{2}+y_{n}^{2}}{2\sigma_{\min}^{4}}\left|\sigma_{n}^{2}(\lambda_{0:n})-\sigma_{n}^{2}(\lambda_{0:n}')\right|
\end{eqnarray*}
where $\sigma_{n}^{2}$ is contracting, in the sense that, for $n\geq
p$,
\begin{eqnarray*}
\left|\sigma_{n}^{2}(\lambda_{0:n})-\sigma_{n}^{2}(\lambda_{0:n}')\right| & = & \left\{ \prod_{i=0}^{p-1}\gamma(\lambda_{n-i})\right\} \left|\sigma_{n-p}^{2}(\lambda_{0:n-p})-\sigma_{n-p}^{2}(\lambda_{0:n-p}')\right|\\
 & \leq & 2\gamma_{\max}^{p}\sigma_{\max}^{2}.
\end{eqnarray*}
Thus, using (\ref{ineq:logdens}), and the fact that $(e^{x}-1)/x$
is an increasing function, Hypothesis 2 is verified with $\tau=\gamma_{\max}$
and \[
\phi_{n}=\tau^{-q}\left[\exp\left\{ \frac{\tau^{q}\sigma_{\max}^{2}\left(\sigma_{\min}^{2}+y_{n}^{2}\right)}{\sigma_{\min}^{4}}\right\} -1\right],\]
for any $q\leq p$. Finally, to compute the expectation with respect
to process $\left(Y_{n}\right)$ of the error bound (\ref{eq:error_bound}),
one may use repetitively the following results:
 \[
\mathbb{E}\left[\exp\left(aY_{n}^{2}\right)|Y_{1:n-1}\right]\leq\left(1-2a\sigma_{\max}^{2}\right)^{-1/2}\]
for $a<1/2\sigma_{\max}^{2}$, using standard calculations and the
fact that $Y_{n}$, conditional on $Y_{1:n-1}$ and $\Lambda_{0:n}=\lambda_{0:n}$
is $\mathcal{N}\left(0,\sigma_{n}^{2}(\lambda_{0:n})\right)$. This
implies in particular that: \[
\mathbb{E}\left[a_{n}b_{n}|Y_{1:n-1}\right]\leq\left(2\frac{\sigma_{\min}^{2}}{\sigma_{\max}^{2}}-1\right)^{-1/2}=c\]
where the constant $c$ is well-defined since $\sigma_{\max}^{2}/\sigma_{\min}^{2}<2$,
then by Jensen inequality,
 \[
\mathbb{E}\left[\frac{1}{a_{n}b_{n}}\bigm\vert Y_{1:n-1}\right]\geq c^{-1},\]
and similarly,
 \[
\mathbb{E}\left[\phi_{n}|Y_{1:n-1}\right]\leq\tau^{-q}\left[\exp\left\{ \tau^{q}\frac{\sigma_{\max}^{2}}{\sigma_{\min}^{2}}\right\} \left(1-2\tau^{q}\frac{\sigma_{\max}^{4}}{\sigma_{\min}^{4}}\right)^{-1/2}-1\right]=\phi\]
where $\phi$ is properly defined for $q$ large enough. Using the
above results recursively on the sum on the RHS of (\ref{eq:error_bound}),
one obtains the same expression as in (\ref{eq:error_cor}) for the
error bound than in Corollary \ref{cor:time-unif} for time-uniform
estimates (with the values of $c$, $\phi$, $\tau$ as defined above),
and concludes similarly. 
\end{proof}
If $\beta$ is allowed to take positive values, stability results may
be obtained under more restrictive conditions. In particular, one may
impose that $\gamma$ is a constant function.
\begin{lem}
For the general mixture GARCH model defined above, the expected particle
error is uniformly stable in time, i.e. there exist constants $C$,
$D$, such that
 \[
\mathbb{E}\left[\mathcal{E}_{n,N}^{p}(Y_{1:n})\right]\leq C\left\{ \log(N)+D\right\} \left(\frac{1}{\sqrt{N}}\right)^{1+3\log c/\log\gamma}\]
provided $\gamma$ is a constant function,
$\gamma(\lambda)=\gamma$,  $\tau c^{3}<1$,  $\vartheta<2$,   where
$\tau=\gamma$, $c=\left(2/\vartheta-1\right)^{-1/2}$, $p$ is given by
(\ref{eq:pfN}), and 
\[
\vartheta=\left(\frac{\alpha_{\max}}{\alpha_{\min}}\vee\frac{\beta_{\max}}{\beta_{\min}}\right).\]
\end{lem}
\begin{proof}
We follow the same lines as above, except that the bounds of the process
$\sigma_{n}^{2}(\lambda_{0:n})$ must be replaced by:
 \begin{eqnarray*}
\sigma_{\min}^{2}(n) & = & \frac{\gamma^{n}}{1-\gamma}\alpha_{\min}+\sum_{k=0}^{n-1}(\alpha_{\min}+\beta_{\min}y_{n-k}^{2})\gamma^{k},\\
\sigma_{\max}^{2}(n) & = & \frac{\gamma^{n}}{1-\gamma}\alpha_{\max}+\sum_{k=0}^{n-1}(\alpha_{\max}+\beta_{\max}y_{n-k}^{2})\gamma^{k},\end{eqnarray*}
which, by construction, are such that
 \[
\frac{\sigma_{\max}^{2}(n)}{\sigma_{\min}^{2}(n)}\leq\vartheta<2.\]
Hence, one has again
\[
\mathbb{E}\left[a_{n}b_{n}|Y_{1:n-1}\right]\leq\left(2\frac{\sigma_{\min}^{2}}{\sigma_{\max}^{2}}-1\right)^{-1/2}=c\]
and the rest of the calculation is identical to those of previous
Lemma, with $\tau=\gamma$.
\end{proof}

\subsection{Mixture Kalman model}

We focus on an univariate linear Gaussian model, i.e. conditional
on Markov process $(\Lambda_{n})$, one has $X_{0}=0$ almost surely,
and, for $n\geq1$, \begin{eqnarray*}
X_{n} & = & h(\Lambda_{n})X_{n-1}+\sqrt{w(\Lambda_{n})}W_{n},\\
Y_{n} & = & X_{n}+\sqrt{v(\Lambda_{n})}V_{n},\end{eqnarray*}
 where the $V_{n}$'s and the $W_{n}$'s are independent $\mathcal{N}(0,1)$
variables, and $h$, $v$, $w$ are real-valued functions. Using the
recursions of the Kalman-Bucy Filter \citep{KalBuc}, one is able
to marginalise out the process $X_{n}$, and compute recursively the
probability density of $Y_{n}$, conditional on $\Lambda_{0:n}=\lambda_{0:n}$,
in the following way:
 \[
\Psi_{n}(\lambda_{0:n})=\frac{1}{\sqrt{2\pi\sigma_{n}^{2}(\lambda_{0:n}})}\exp\left[-\frac{\left\{
      y_{n}-\mu_{n}(\lambda_{0:n})\right\}
    ^{2}}{2\sigma_{n}^{2}(\lambda_{0:n})}\right]\]
 where, the following quantities are defined recursively: for $n\geq1$,
\begin{eqnarray}
\mu_{n}(\lambda_{0:n}) & = & h(\lambda_{n})m_{n-1}(\lambda_{0:n-1})\label{eq:mun}\\
\sigma_{n}^{2}(\lambda_{0:n}) & = & h(\lambda_{n})^{2}c_{n-1}(\lambda_{0:n-1})+v(\lambda_{n})+w(\lambda_{n})\label{eq:sign}\\
a_{n}(\lambda_{0:n}) & = & \left\{ h(\lambda_{n})^{2}c_{n-1}(\lambda_{0:n-1})+w(\lambda_{n})\right\} /\sigma_{n}^{2}(\lambda_{0:n})\label{eq:an}\\
m_{n}(\lambda_{0:n}) & = & h(\lambda_{n})m_{n-1}(\lambda_{0:n-1})+a_{n}(\lambda_{0:n})\left\{ y_{n}-\mu_{n}(\lambda_{0:n})\right\} \label{eq:mn}\\
c_{n}(\lambda_{0:n}) & = & h(\lambda_{n})^{2}c_{n-1}(\lambda_{0:n-1})+w(\lambda_{n})-a_{n}(\lambda_{0:n})^{2}\sigma_{n}^{2}(\lambda_{0:n})\label{eq:cn}
\end{eqnarray}
and $m_{0}(\lambda_{0})=c_{0}(\lambda_{0})=0$. 

We make the following assumptions: 
\begin{enumerate}
\item Functions $v$ and $w$ are bounded as follows: for all
  $\lambda\in E$,
\[
0<\underline{v}\leq v(\lambda)\leq\bar{v},\qquad0<\underline{w}\leq w(\lambda)\leq\bar{w}.\]

\item Function $h$ is bounded as follows: for all $\lambda\in E$,
\[
\left|h(\lambda)\right|\leq\overline{h}<1\]

\end{enumerate}
We first prove  the following intermediate results.
\begin{lem}
\label{lem:sig}The sequence $ $$\sigma_{n}^{2}$ is bounded and
uniformly contracting, i.e. for all $p\geq1,$ for all $\lambda_{0:n}$,
$\lambda'_{0:n}$, such that $\lambda_{n-p+1:n}=\lambda_{n-p+1:n}'$,
one has \[
\underline{\sigma}^{2}\leq\sigma_{n}^{2}(\lambda_{0:n})\leq\bar{\sigma}^{2}\quad\left|\sigma_{n}^{2}(\lambda_{0:n})-\sigma_{n}^{2}(\lambda_{0:n}')\right|\leq C_{\sigma}\tau_{\sigma}^{p}\]
where $ $$\underline{\sigma}^{2}=\underline{v}+\underline{w}$, $\bar{\sigma}^{2}=(\bar{h}^{2}+1)\bar{v}+\bar{w}$,
$C_{\sigma}=\bar{h}^{2}\bar{v}/\tau_{\sigma}$, and
 \[
\tau_{\sigma}=\frac{1}{1+\underline{w}/\bar{v}+2\sqrt{\underline{w}/\bar{v}+\underline{w}^{2}/\bar{v}^{2}}}<1.\]
\end{lem}
\begin{proof}
From (\ref{eq:cn}), one deduces that
\begin{equation}
\frac{1}{c_{n}(\lambda_{0:n})}=\frac{1}{v(\lambda_{n})}+\frac{1}{h(\lambda_{n})^{2}c_{n-1}(\lambda_{0:n-1})+w(\lambda_{n})}\label{eq:invcn}\end{equation}
 thus \[
\left(\frac{1}{\overline{v}}+\frac{1}{\overline{w}}\right)^{-1}\leq c_{n}(\lambda_{0:n})\leq\bar{v}\]
and, from (\ref{eq:sign}), that
 \[
\underline{v}+\underline{w}\leq\sigma_{n}^{2}(\lambda_{0:n})\leq(\bar{h}^{2}+1)\bar{v}+\bar{w}.\]
 In addition, (\ref{eq:invcn}) implies that
 \[
\log\left\{ c_{n}(\lambda_{0:n})\right\} =\Upsilon\left(\log\left\{ c_{n-1}(\lambda_{0:n-1})\right\} ,\lambda_{n}\right)\]
 where \[
\Upsilon\left(c,\lambda\right)=-\log\left\{ \frac{1}{v(\lambda)}+\frac{1}{h(\lambda)^{2}e^{c}+w(\lambda)}\right\} .\]
 It is easy to show that, for a fixed $\lambda$, the derivative of
$\Upsilon\left(c,\lambda\right)$ with respect to $c$ is bounded
from above by $\tau_{\sigma}$ as defined above. Thus, $\Upsilon\left(c,\lambda\right)$
is a contracting function, and, by induction, for $n\geq p$,
\begin{eqnarray*}
\left|\sigma_{n}^{2}(\lambda_{0:n})-\sigma_{n}^{2}(\lambda_{0:n}')\right| 
& = & \left|h(\lambda_{n})\right|^2\left|c_{n-1}(\lambda_{0:n-1})-c_{n-1}(\lambda_{0:n-1}')\right|\\
 & \leq & \bar{h}^2\bar{v}\left|\log c_{n-1}(\lambda_{0:n-1})-\log c_{n-1}(\lambda_{0:n-1}')\right|\\
 & \leq & C_{\sigma}\tau_{\sigma}^{p}.
\end{eqnarray*}
where $\tau_{\sigma}$ and $C_{\sigma}$ were defined above.\end{proof}
\begin{lem}
The sequence $ $$\mu_{n}$ is bounded and contracting in the sense
that there exists $C_{\mu}>0$ such that, for all $p\geq1,$ for all
$n\geq p$, and $\lambda_{0:n}$, $\lambda'_{0:n}$, such that $\lambda_{n-p+1:n}=\lambda'_{n-p+1:n}$,
one has
 \[
\left|\mu_{n}(\lambda_{0:n})\right|\leq\frac{\bar{a}\bar{h}}{1-\tilde{a}\bar{h}}M_{n-1},
\quad\left|\mu_{n}(\lambda_{0:n})-\mu_{n}(\lambda_{0:n}')\right|\leq
C_{\mu}M_{n-1}\tau^{p},
\]
where
 \[ M_{n}=\max_{i=1,\ldots n}\left(\left|y_{i}\right|\right),\quad
\tau=\tau_{\sigma}\vee\bar{h},\quad\bar{a}=\left(1-\frac{\underline{v}}{\bar{h}^{2}\bar{v}+\bar{w}+\underline{v}}\right),\quad\tilde{a}=\frac{\bar{v}}{\bar{v}+\underline{w}}.\]
\end{lem}
\begin{proof}
Note first that \[
1-\tilde{a}=\frac{\underline{w}}{\bar{v}+\underline{w}}\leq a_{n}(\lambda_{0:n})\leq\left(1-\frac{\underline{v}}{\bar{h}^{2}\bar{v}+\bar{w}+\underline{v}}\right)=\bar{a}\]
so one shows recursively, using (\ref{eq:mun}) and (\ref{eq:mn}),
that:\[
\left|\mu_{n}(\lambda_{0:n})\right|\leq\frac{\bar{a}\bar{h}}{1-\tilde{a}\bar{h}}M_{n-1}\]
and that, for $\lambda_{0:n}$, $\lambda'_{0:n}$ such that $\lambda_{n-p+1:n}=\lambda'_{n-p+1:n}$,
\begin{eqnarray}
\lefteqn{\left|\mu_{n}(\lambda_{0:n})-\mu_{n}(\lambda_{0:n}')\right|}\label{eq:diffmun}\\
 & \leq &
 M_{n-1}\left[\sum_{i=1}^{p}\bar{h}^{i}\left|a_{n-i}\prod_{j=1}^{i-1}(1-a_{n-j})
-a_{n-i}'\prod_{j=1}^{i-1}(1-a_{n-j}')\right|+2\bar{h}^{p+1}\right] \nonumber
\end{eqnarray}
where $a_{n-i}$, $a_{n-i}'$ are short-hands for $a_{n-i}(\lambda_{0:n-i})$,
$a_{n-i}(\lambda_{0:n-i}')$. The sequence $a_{n}$ itself in contracting,
since, from (\ref{eq:an}), one has, for $i<p$:
\begin{eqnarray*}
  \left|a_{n-i}-a_{n-i}'\right| 
  & \leq & \frac{\bar{v}}{\underline{\sigma}^{4}}\left|\sigma_{n-i}^{2}(\lambda_{0:n-i})-\sigma_{n-i-1}^{2}(\lambda_{0:n-i}')\right|\\
  & \leq &
  \frac{\bar{v}C_{\sigma}}{\underline{\sigma}^{4}}\tau_{\sigma}^{p-i}
\end{eqnarray*}
so (\ref{eq:diffmun}) and the fact that
 $\left|xy-x'y'\right|\leq|x-x'|+|y-y'|$
provided $x,x',y,y'\in[0,1]$ leads to
 \begin{eqnarray*}
\lefteqn{\left|\mu_{n}(\lambda_{0:n})-\mu_{n}(\lambda_{0:n}')\right|}\\
 & \leq & M_{n-1}\left[\frac{\bar{v}C_{\sigma}}{\underline{\sigma}^{4}}\sum_{i=1}^{p}\bar{h}^{i}\left(\tau_{\sigma}^{p-i}+\ldots+\tau_{\sigma}^{p-1}\right)+2\bar{h}^{p+1}\right]\\
 & \leq & M_{n-1}\left[\frac{\bar{v}C_{\sigma}\tau_{\sigma}^{p-1}}{\underline{\sigma}^{4}}\sum_{i=1}^{p}\bar{h}^{i}\left(\frac{\tau_{\sigma}^{-i}-1}{\tau_{\sigma}^{-1}-1}\right)+2\bar{h}^{p+1}\right]\\
 & \leq & M_{n-1}C_{\mu}\tau^{p}\end{eqnarray*}
for $\tau=\tau_{\sigma}\vee\bar{h}$, and a well chosen value of $C_{\mu}$. 
\end{proof}
We are now able to state the main result. 
\begin{lem}
For the model above, the particle error is bounded uniformly in time,
i.e. there exist $C$, $D$, such that
 \[
\mathcal{E}_{n,N}^{p}(y_{1:n})\leq C\left\{ \log(N)+D\right\} \left(\frac{1}{\sqrt{N}}\right)^{1+3\log c/\log\tau},\]
almost surely, for $p$ given by (\ref{eq:pfN}), provided the
realizations $y_{n}$ are bounded, i.e. $\left|y_{n}\right|\leq C_{y}$
for all $n\geq1$, and that $\tau c^{3}<1,$ with $\tau=\bar{h}\vee\tau_{\sigma}$
and \[
c=\frac{\bar{\sigma}}{\underline{\sigma}}\exp\left[\frac{C_{y}^{2}}{\underline{\sigma}^{2}}\left\{
    1+\left(\frac{\bar{a}\bar{h}}{1-\tilde{a}\bar{h}}\right)^{2}\right\}
\right],
\quad \tau_{\sigma}=\frac{1}{1+\underline{w}/\bar{v}+2\sqrt{\underline{w}/\bar{v}+\underline{w}^{2}/\bar{v}^{2}}}<1.\]
\end{lem}
\begin{proof}
This proposition is a direct application of Corrolary \ref{cor:time-unif},
so we need only to prove that Hypotheses 2 and 3 are fulfilled. For
Hypothesis 2, one may take
\[
\frac{1}{a_{n}}=\frac{1}{\sqrt{2\pi\bar{\sigma}^{2}}}\exp\left[-\frac{C_{y}^{2}}{\underline{\sigma}^{2}}\left\{ 1+\left(\frac{\bar{a}\bar{h}}{1-\tilde{a}\bar{h}}\right)^{2}\right\} \right],\quad b_{n}=\frac{1}{\sqrt{2\pi\underline{\sigma}^{2}}}\]
so that $a_{n}b_{n}\leq c$ for $c$ defined above. For Hypothesis
3, one has:
 \begin{eqnarray*}
2\left|\log\Psi_{n}(\lambda_{0:n})-\log\Psi_{n}(\lambda_{0:n}')\right| & \leq & \left|\log\sigma_{n}^{2}(\lambda_{0:n})-\log\sigma_{n}^{2}(\lambda_{0:n}')\right|\\
 &  & +\left|\frac{\left\{ y_{n}-\mu_{n}(\lambda_{0:n})\right\} ^{2}}{\sigma_{n}^{2}(\lambda_{0:n})}-\frac{\left\{ y_{n}-\mu_{n}(\lambda_{0:n}')\right\} ^{2}}{\sigma_{n}^{2}(\lambda_{0:n}')}\right|\end{eqnarray*}
where the first term is such that
 \begin{eqnarray*}
\left|\log\sigma_{n}^{2}(\lambda_{0:n})-\log\sigma_{n}^{2}(\lambda_{0:n}')\right| & \leq & \frac{1}{\underline{\sigma}^{2}}\left|\sigma_{n}^{2}(\lambda_{0:n})-\sigma_{n}^{2}(\lambda_{0:n}')\right|\\
 & \leq & \frac{C_{\sigma}}{\underline{\sigma}^{2}}\tau_{\sigma}^{p}\end{eqnarray*}
according to Lemma \ref{lem:sig}, and the second term is such that
 \begin{eqnarray*}
\lefteqn{\left|\frac{\left\{ y_{n}-\mu_{n}(\lambda_{0:n})\right\} ^{2}}{\sigma_{n}^{2}(\lambda_{0:n})}-\frac{\left\{ y_{n}-\mu_{n}(\lambda_{0:n}')\right\} ^{2}}{\sigma_{n}^{2}(\lambda_{0:n}')}\right|}\\
 & \leq & \frac{1}{\sigma_{n}^{2}(\lambda_{0:n}')}\left|\left\{ y_{n}-\mu_{n}(\lambda_{0:n})\right\} ^{2}-\left\{ y_{n}-\mu_{n}(\lambda_{0:n}')\right\} ^{2}\right|\\
 &  & +\frac{\left\{ y_{n}-\mu_{n}(\lambda_{0:n})\right\}
   ^{2}}{\sigma_{n}^{2}(\lambda_{0:n})\sigma_{n}^{2}(\lambda_{0:n}')}
\left\vert\sigma_{n}^{2}(\lambda_{0:n})-\sigma_{n}^{2}(\lambda_{0:n}')\right\vert\\
 & \leq & \frac{2C_{\mu}C_{y}^{2}}{\underline{\sigma}^{2}}\left(1+\frac{\bar{a}\bar{h}}{1-\tilde{a}\bar{h}}\right)\tau^{p}+\frac{2C_{y}^{2}C_{\sigma}}{\underline{\sigma}^{4}}\left[1+\left(\frac{\bar{a}\bar{h}}{1-\tilde{a}\bar{h}}\right)^{2}\right]\tau_{\sigma}^{p}\end{eqnarray*}
 and one concludes using (\ref{ineq:logdens}) and taking
 \[
\phi=\phi_{n}=\exp\left\{ \frac{C_{\sigma}}{2\underline{\sigma}^{2}}+\frac{C_{\mu}C_{y}^{2}}{\underline{\sigma}^{2}}\left(1+\frac{\bar{a}\bar{h}}{1-\tilde{a}\bar{h}}\right)+\frac{C_{y}^{2}C_{\sigma}}{\underline{\sigma}^{4}}\left[1+\left(\frac{\bar{a}\bar{h}}{1-\tilde{a}\bar{h}}\right)^{2}\right]\right\} -1.\]

\end{proof}
Obviously, the boundness condition on the realizations $y_{n}$ is
not entirely satisfactory, as the generating process of $(Y_{n})$
is such that $Y_{n}$ should leave any interval eventually. However,
$Y_{n}$ is marginally a Gaussian variable with variance uniformly
bounded in time (since $\bar{h}<1$), so this remains a reasonable
approximation if $C_{y}$ is large enough. Generalizing the above
result to more general conditions is left for future research.

\subsection{Application to standard state-space models}

Consider a `standard' state-space model, based on a linear auto-regressive
state process $(X_{n})$: \begin{equation}
X_{n}=\rho X_{n-1}+\Lambda_{n},\quad\Lambda_{1},\ldots,\Lambda_{n},\ldots\:\mbox{i.i.d.}\label{eq:autoreg}\end{equation}
for $t\geq0$, $\rho\in(-1,1)$ and $X_{0}=\Lambda_{0}$, and an observed
process $(Y_{n})$, with conditional density, with respect to an appropriate
dominating measure, and conditional on $X_{n}=x_{n}$, given by the
potential function $\Psi_{n}^{X}(x_{n})$. 

In this section, we show how to apply our stability results to such
a standard state-space model, where the potential function depends
only on the current state $X_{n}$. We rewrite the model as a state
space model with hidden Markov chain $(\Lambda_{n})$, and observed
process $(Y_{n})$ corresponding to potential function
 \[
\Psi_{n}(\lambda_{0:n})=\Psi_{n}^{X}\left(\sum_{k=0}^{n}\rho^{k}\lambda_{n-k}\right),\]
where the argument $x_{n}$ in the right hand side has been substituted
with the appropriate function of $\lambda_{0:n}$, as derived from
(\ref{eq:autoreg}).

Clearly, the reformulated model satisfies Hypothesis 1: the $(\Lambda_{n})$
are i.i.d., hence they form a Markov chain with mixing coefficient
$\varepsilon_{n}=1$. If we assume that the $\Psi_{n}(\lambda_{0:n})$
are such that Hypotheses 2 and 3 hold as well, then we can apply directly
Theorem \ref{prop:maj}. However, the path-dependent formulation
of this model is artificial, and, in practice, we are interested in
filtering the process $X_{n}$, conditional on the $Y_{n}$'s, rather
than filtering the $\Lambda_{n}$'s, again conditional on the $Y_{n}$'s.
More precisely, we wish to approximate the conditional expectation
of\[
g(X_{n})=g\left(\sum_{k=0}^{p-1}\rho^{k}\lambda_{n-k}+\sum_{k=p}^{n}\rho^{k}\lambda_{n-k}\right),\]
for some bounded function $g$, and, provided $g$ is also Lipschitz,
with constant $K$, and that the $\lambda_{n}$'s  lie in interval
$[-l,l],$ for some $l\geq0$, one has:
\[
\left|g\left(\sum_{k=0}^{p-1}\rho^{k}\lambda_{n-k}+\sum_{k=p}^{n}\rho^{k}\lambda_{n-k}\right)-g\left(\sum_{k=0}^{p-1}\rho^{k}\lambda_{n-k}\right)\right|\leq\frac{Kl}{1-\tau}\tau^{p},\]
where $\tau=\left|\rho\right|$. Therefore, we must consider an additional
term in the particle error attached to the filtering of $(X_{n})$,
which stems from the difference between the filtering distribution
of $X_{n}$ and that of $\Lambda_{n-p+1:n}$, for some integer $p$.
Consider the following estimate of the particle error for functions
of $X_{n}$:
\[
\mathcal{E}_{n,N}^{X}(y_{1:n})=\sup_{g:\Vert g\Vert_{\infty}=1,g\in Lip(K)}\mathbb{E}_{N}\left(\left|\langle R_{1:n}\zeta-R_{1:n}^{N}\zeta,f_{g}\rangle\right|\ |Y_{0:n}=y_{0:n}\right)\]
where $Lip(K)$ denotes the set of Lipschitz functions with Lipschitz
constant $K$, and $f_{g}$ is the function $E^{n+1}\rightarrow\mathbb{R}$
such that \[
f_{g}(\lambda_{0:n})=g\left(\sum_{k=0}^{n}\rho^{k}\lambda_{n-k}\right),\]
i.e., loosely speaking, $f_{g}(\lambda_{0:n})=g(x_{n})$, where $x_n$
must be substituted by its expression as a function of $\lambda_{0:n}$. 
\begin{lem}
For the state-space model described above, one has, for any $n\geq p$,
 \[
\mathcal{E}_{n,N}^{X}(y_{1:n})\leq\mathcal{E}_{n,N}^{p}(y_{1:n})+\frac{Kl}{1-\tau}\tau^{p}.\]

\end{lem}
Taking into account this additional error term, we can derive time-uniform
estimates of the stability of the particle algorithm. For the sake
of space, we focus on the following simple example: $Y_{n}\in\{-1,1\},$
$Y_{n}=1$ with probability $1/(1+e^{X_{n}})$, $Y_{n}=-1$ otherwise.
The potential function (for the model in its standard formulation)
equals:\[
\Psi_{n}^{X}(x_{n})=\frac{1}{1+e^{y_{n}x_{n}}}.\]

We recall that the support of the $(\Lambda_{n})$ is $[-l,l]$, and
therefore $X_{n}\in[-l',l']$ almost surely, with $l'=l/(1-\tau)$.
Thus,  Hypothesis 3 holds for $b_{n}=1/(1+e^{-l'})$, $a_{n}=1+e^{l'}$.
For Hypothesis 2, standard calculations show that, for two vectors
$\lambda_{0:n}$ and $\lambda_{0:n}'$ such that $\lambda_{n-p+1:n}=\lambda_{n-p+1:n}'$,
one has \begin{eqnarray*}
\left|\log\Psi_{n}(\lambda_{0:n})-\log\Psi_{n}(\lambda_{0:n}')\right| 
& \leq & \left|\sum_{k=p}^{n}\rho^{k}(\lambda_{n-k}-\lambda'_{n-k})\right|\\
 & \leq & 2l'\tau^{p}\end{eqnarray*}
provided $\tau=\left|\rho\right|$. Hence, using (\ref{ineq:logdens})
inequality, Hypothesis 2 holds, with $\phi_{n}=e^{2l'}-1.$ 

For this specific model,  we have the following result. 
\begin{lem}
For the specific model described above, and provided $c\tau^{3}<1$,
where $\tau=\left|\rho\right|$, $c=e^{l'}$, one
has: \[
\mathcal{E}_{n,N}^{X}(y_{1:n})\leq C\left\{ \log(N)+D\right\} \left(\frac{1}{\sqrt{N}}\right)^{1+3\log c/\log\tau}+\frac{E}{\sqrt{N}}\]
where $C$ and $D$ were defined in Corollary \ref{cor:time-unif},
$\phi=e^{2l'}-1$, and $E=4Kl'c/3\phi$. 
\end{lem}
The above model does not fulfil the usual conditions required in standard
stability results, see e.g. \citet[Section 7.4.3]{delMoral:book},
because the Markov chain $(X_{n})$ is not mixing. Thus, it is remarkable
that the time-uniform stability of this model is established using
a Feynman-Kac formulation with path-dependent potentials.

\section{Conclusion}

To extend our results to a broader class of models, three directions
may be worth investigating. First, it may be possible to bound directly
the particle error, without resorting to a comparison with an artificial,
truncated potential function. It seems difficult however to avoid
some form of truncation, as the path process $\Lambda_{0:n}$ itself
does not benefit to any sort of mixing property, while fixed segments
$\Lambda_{n-p+1:n}$ do. Second, one may try to loosen Hypothesis
1 (Markov kernel is mixing) and Hypothesis 3 (potential function is
bounded), using for instance \citet{OudRub}'s approach. Third, it
seems possible to adapt our general result on the particle error bound
to several models not considered in this paper, in particular standard
models with potential functions depending on the last state only, by 
using and extending the approach developed in the previous Section. 

\section*{Acknowledgements}

Part of this work was performed while the authors were invited
by the SAMSI institute during the 2008-2009 Program on Sequential
Monte Carlo Methods.

\bibliographystyle{apalike}
\bibliography{complete}

\begin{thebibliography}{}

\bibitem[Andrieu and Doucet, 2002]{AndDou}
Andrieu, C. and Doucet, A. (2002).
\newblock Particle filtering for partially observed {G}aussian state space
  models.
\newblock {\em J. R. Statist. Soc. B}, 64(4):827--836.

\bibitem[Atar and Zeitouni, 1997]{atar1997esn}
Atar, R. and Zeitouni, O. (1997).
\newblock {Exponential stability for nonlinear filtering}.
\newblock 33(6):697--725.

\bibitem[Bollerslev, 1986]{bollerslev1986gac}
Bollerslev, T. (1986).
\newblock {Generalized autoregressive conditional heteroskedasticity}.
\newblock {\em Journal of Econometrics}, 31(3):307--327.

\bibitem[Capp\'e et~al., 2005]{CapMouRyd}
Capp\'e, O., Moulines, E., and Ryd\'en, T. (2005).
\newblock {\em Inference in Hidden {M}arkov Models}.
\newblock Springer-Verlag, New York.

\bibitem[Chen and Liu, 2000]{ChenLiu}
Chen, R. and Liu, J. (2000).
\newblock Mixture {K}alman filters.
\newblock {\em J. R. Statist. Soc. B}, 62:493--508.

\bibitem[Chopin, 2007]{Chopin:change}
Chopin, N. (2007).
\newblock Dynamic detection of change points in long time series.
\newblock {\em Ann. of the Inst. of Stat. Math.}, 59(2):349--366.

\bibitem[Del~Moral, 2004]{delMoral:book}
Del~Moral, P. (2004).
\newblock {\em Feynman-Kac Formulae}.
\newblock Springer.

\bibitem[Doucet et~al., 2001]{DouFreiGor}
Doucet, A., de~Freitas, N., and Gordon, N.~J. (2001).
\newblock {\em Sequential {M}onte {C}arlo Methods in Practice}.
\newblock Springer-Verlag, New York.

\bibitem[Doucet et~al., 2000]{DouGodAnd}
Doucet, A., Godsill, S., and Andrieu, C. (2000).
\newblock On sequential {M}onte {C}arlo sampling methods for {B}ayesian
  filtering.
\newblock {\em Statist. Comput.}, 10(3):197--208.

\bibitem[Kalman and Bucy, 1961]{KalBuc}
Kalman, R.~E. and Bucy, R.~S. (1961).
\newblock New results in linear filtering and prediction theory.
\newblock {\em Trans. Amer. Soc. Mech. Eng., J. Basic Eng.}, 83:95--108.

\bibitem[K{\"u}nsch, 2001]{Kun:SSHMM}
K{\"u}nsch, H. (2001).
\newblock State space and hidden {M}arkov models.
\newblock In Barndorff-Nielsen, O.~E., Cox, D.~R., and Kl{\"u}ppelberg, C.,
  editors, {\em Complex Stochastic Systems}, pages 109--173. Chapman and Hall.

\bibitem[Le~Gland and Oudjane, 2004]{legland2004sau}
Le~Gland, F. and Oudjane, N. (2004).
\newblock {Stability and uniform approximation of nonlinear filters using the
  Hilbert metric and application to particle filters}.
\newblock {\em Ann. Appl. Prob.}, 14(1):144--187.

\bibitem[Oudjane and Rubenthaler, 2005]{OudRub}
Oudjane, N. and Rubenthaler, S. (2005).
\newblock Stability and uniform particle approximation of nonlinear filters in
  case of non ergodic signals.
\newblock {\em Stochastic Analysis and applications}, 23:421--448.

\end{thebibliography}

\end{document}